%% file: ComplexMult.tex
\documentclass[11pt]{amsart}

\input{MTemplate}

\usepackage{mathrsfs}

\title{Complex Multiplication tests for Elliptic Curves}
\author{Denis Xavier Charles}
\thanks{Research supported in part by NSF grant CCR-9988202.}
\address{Department of Computer Science, University of Wisconsin-Madison, Madison WI
 - 53706.}

\email{cdx@cs.wisc.edu}

\date{17 May, 2004}

\begin{document}
\begin{abstract}
We consider the problem of checking whether an elliptic curve defined over a given number field
has complex multiplication. We study two polynomial time algorithms for this problem, one randomized
and the other deterministic. The randomized algorithm can be adapted to yield 
the discriminant of the endomorphism ring of the curve.\\

\vspace{0.5cm}
\parindent=0pt
{\sc Keywords. } Algorithms, Elliptic Curves, Complex Multiplication, Endomorphism Ring, $\ell$-adic representations,
Chebotarev Density Theorem.
\end{abstract}
\maketitle

\section{Introduction}

It is a well known fact that the endomorphism ring of an elliptic curve
over a {\em number field} is isomorphic to either $\mathbb{Z}$
or an order in an imaginary quadratic field. If the latter holds then
the curve is said to have {\em complex multiplication} (CM.) Elliptic curves with
complex multiplication have found applications in cryptography and coding theory, 
since there are closed form expressions for the number of points on such curves
modulo prime ideals. 
This property was also utilized in the Atkin-Morain
primality proving method \cite{atmor93}. Constructing elliptic curves
with complex multiplication is computationally very expensive.
In this article we show that testing an elliptic curve for CM is easy.\\

If one fixes the number field over which the curves are defined, then CM testing
becomes very easy, albeit with considerable pre-computation.
For this reason we consider the number field as being part of the input (this issue is
explained in section \S\ref{sec_prob:this}). 
Once one defines the problem in this way, an approach
immediately suggests itself: transform the method of constructing curves with complex
multiplication into a solution for this problem. Unfortunately, to implement
this method one needs {\em good} effective lower bounds on class numbers of imaginary quadratic
fields, which is a notorious open problem. This approach and its analysis is the subject of
\S\ref{sec_direct:this}. \\

Our next approach, discussed in \S\ref{sec_rand:this}, 
uses the elegant results of Deuring on the reduction of endomorphism
rings of elliptic curves and Serre on the density of supersingular
primes. The approach is based on the observation that supersingular primes
are plentiful for curves with complex multiplication. This yields a two-sided
error probabilistic polynomial time algorithm for this problem. We also show how this
method can be adapted to find the discriminant of the endomorphism ring, but the analysis of
this stage of the algorithm presents some challenging open questions. However, we can use
the results we obtain here to make the error in the randomized algorithm one-sided. A similar algorithm
is sketched in \cite{cnst98} without a precise analysis of the probability of failure
and the running time. We improve their results in two ways. First, our algorithm is
simpler to implement. Second, unlike theirs, our proof is rigorous and does not
rely on unproven heuristic assumptions.\\

The final method, which we believe is new, 
discussed in \S\ref{sec_determ:this} is based on studying the image of the galois representations
afforded by $\ell$-torsion points on the curve. This method is deterministic and has a polynomial running time,
but we are unable to bound the (multiplicative) constant in the running time effectively.

\section{Preliminaries}\label{sec_prelim:this}
Let $L$ be a number field and let $E/L$ be an elliptic curve. Every elliptic curve over $L$
is isomorphic over $L$ to one that is given by an equation of the form (\cite{sil86} III.\S1)
\begin{align}\label{eqn_wform:this}
	Y^2 Z = X^3 + AXZ^2 + BZ^3
\end{align}
with $A,B \in L$ and $4A^3 + 27B^2 \neq 0$. If $E$ is an elliptic curve that is given by an
equation of the above form, then we define the {\bf discriminant} of $E$
by
\begin{align*}
	\Delta_E = -16(4A^3 + 27B^2)
\end{align*}
and the {\bf j-invariant} of $E$ to be the quantity
\begin{align*}
	j_E = \frac{-1728 (4A)^3}{\Delta_E}.
\end{align*}

For the rest of the article, 
an elliptic curve over a number field $L$ is a curve given by
an equation of the form (\ref{eqn_wform:this}) with coefficients in $L$.

\subsection{Structure of the Endomorphism ring}
Let $E_1,E_2$ be two elliptic curves defined over $L$. $\mathrm{Hom}(E_1,E_2)$
is the set $\{\phi ~|~ \phi : E_1 \rightarrow E_2 \text{ is an isogeny}\}$.
$\mathrm{Hom}(E_1,E_2)$ is given a group structure by defining addition of maps pointwise.
$\mathrm{End}(E)$ as a set is defined to be $\mathrm{Hom}(E,E)$. $\mathrm{End}(E)$
is a {\bf ring} with multiplication defined to be composition of isogenies.
The multiplication-by-$m$ map $[m]$ belongs to $\mathrm{End}(E)$ for each $m \in \mathbb{Z}.$ 
In fact, the map $\mathbb{Z} \rightarrow \mathrm{End}(E)$
given by $m \mapsto [m]$ is an injection of rings. The following result of Deuring gives the
possibilities for $\mathrm{End}(E)$.

\begin{Thm}[Deuring] Let $E/L$ be an elliptic curve, then $\mathrm{End}(E)$ is either $\mathbb{Z}$
or $\mathscr{O}$, an order in an imaginary quadratic field $K$.
\end{Thm}

Suppose $E/L$ is an elliptic curve with $\mathscr{O} = \mathrm{End}(E) \neq \mathbb{Z}$. Then
we say that $E$ has {\bf complex multiplication} (by $\mathscr{O}$.) Sometimes, for brevity,
we write ``E has CM'' instead of ``E has complex multiplication.''

\subsection{Weil Height} We introduce the notion of the Weil height of an algebraic number
which we need in \S\ref{sec_determ:this}.

\begin{Def} Let $\alpha \in \overline{\mathbb{Q}}$ be an algebraic number with minimal polynomial
\begin{align*}
p_{\alpha}(x) = a_0 x^d + a_1 x^{d-1} + \cdots + a_d \in \mathbb{Z}[x].
\end{align*}
Assume that $p_{\alpha}(x) = a_0(x-\alpha_1)(x-\alpha_2)\cdots(x-\alpha_d)$ with $\alpha_i \in \mathbb{C}$.
Then the {\bf absolute logarithmic Weil height} (or just Weil height) of $\alpha$ is defined
to be the quantity
\begin{align*}
\mathbf{h}(\alpha) = \frac{1}{d}\biggl( \log |a_0| + \sum_{1 \leq i \leq d} \max\{1,|\alpha_i|\}\biggr).
\end{align*}
\end{Def}

With the notation of the definition, we have the following useful bound (\cite{fel82} Lemma 8.2)
\begin{align*}
	\mathbf{h}(\alpha) \leq \frac{1}{d} \log \sum_{i} |a_i|.
\end{align*}

Thus the Weil height of an algebraic number is bounded polynomially by the encoding
length of its minimal polynomial. Also, we denote the quantity $\sum_i |a_i|$
by $\mathbf{w}(\alpha)$.\\

If $E/L$ is an elliptic curve we define the {\bf Weil height} of $E$ to be $\mathbf{h}(j_E)$,
the Weil height of its $j$-invariant.

\section{The Problem}\label{sec_prob:this}
The computational problem that is the focus of this article is the following:\\

{\bf Complex multiplication of elliptic curves:}\\

{\bf Input: } A number field $L$, and an elliptic curve $E: Y^2Z = X^3 + AXZ^2 + BZ^3$
with $A,B \in L$.\\
{\bf Question: } Does $E$ have complex multiplication?\\

We will assume that $L = \mathbb{Q}(j_E)$, since $E$ always has a model over $\mathbb{Q}(j_E)$
and we can restrict to the subfield generated by $j_E$. The input is specified by
giving the minimal polynomial of $A$ and $B$ from which the minimal polynomial of $j_E$
can be determined efficiently. The size of the input is measured by the size of the encoding
of the minimal polynomials of $A$ and $B$. The encoding length of a polynomial
$p(x) = a_0 x^d + a_1 x^{d-1} + \cdots + a_d$, with integer coefficients, is defined to be the quantity 
$\sum_{0 \leq i \leq d} \max\{1,\log |a_i|\}$. Note that the encoding length of a non-zero polynomial $p(x)$
is at least the degree of $p(x).$ \\

Our main concern is the complexity of the above decision problem. A consequence of the
algorithms presented in this article is that the above decision problem is in ${\sf P}$.
Next, we explain why the number field needs to be part of the input.\\

The complex points on $E$, namely $E(\mathbb{C})$, has a particularly
simple interpretation as $\mathbb{C}/\mathcal{L}_E$, where $\mathcal{L}_E$
is a rank $2$ lattice such that $\mathcal{L}_E \otimes_{\mathbb{Z}}\mathbb{R} = \mathbb{C}$.
In this description, isomorphic elliptic curves correspond to lattices that differ by a non-zero
complex scalar (\cite{sil86} VI Ex. 6.6). Suppose $E/\mathbb{C}$ is
given by a lattice $\mathcal{L}_E$, then there is an isomorphic elliptic curve given by the
lattice $\mathbb{Z}+\mathbb{Z}\tau_E$ with $\tau_E \in \mathfrak{H}$, where 
$\mathfrak{H} = \{z \in \mathbb{C} ~:~ \Im z > 0\}.$ There is a simple criterion
for deciding when $E$ has complex multiplication, provided $E$ is given as $\mathbb{C}/(\mathbb{Z}+\mathbb{Z}\tau_E)$
(\cite{sil86} Theorem VI.5.5):\\

Let $\tau$ be an imaginary quadratic number with minimal polynomial $ax^2 + bx + c$ and 
$\gcd(a,b,c) = 1.$ Then the {\bf discriminant} of $\tau$ is $b^2 - 4ac$.

\begin{Thm}\label{thm_cm_crit:this}
Let $E \cong \mathbb{C}/(\mathbb{Z} + \mathbb{Z}\tau_E)$ with $\tau_E \in \mathfrak{H}$. 
Then $E$ has complex multiplication by an order $\mathscr{O}_D$ of discriminant $D$
iff $\tau_E$ is a quadratic number of discriminant $D$ as defined above.
\end{Thm}

We also have the following important theorem (see \cite{coh93} Theorem 7.2.14 or \cite{sil94} Chapter 2):

\begin{Thm} \label{thm_hilb:this}
Let $\tau \in \mathfrak{H}$ be an imaginary quadratic number, and let $D$ be its discriminant.
Then $j(\tau)$ (here $j$ is the usual modular $j$-function) is an algebraic integer of degree equal to $h(D)$,
where $h(D)$ is the class number of the imaginary quadratic order of discriminant $D$. More precisely,
the minimal polynomial of $j(\tau)$ over $\mathbb{Z}$ is the equation $\prod(X - j(\alpha))$,
where $\alpha$ runs over the quadratic numbers associated to the reduced forms of discriminant $D$.
\end{Thm}

We can interpret Theorems \ref{thm_cm_crit:this} and \ref{thm_hilb:this} as follows.
If $E/L$ has complex multiplication by $\mathscr{O}_D$, an order of discriminant $D$,
then its $j$-invariant has only $h(D)$ possibilities, and is an
algebraic integer of degree $h(D)$. Noting that $h(D) \rightarrow \infty$
as $D \rightarrow -\infty$, one concludes that if we fix a number field $L$,
then there are only finitely many $j$-invariants of elliptic curves defined over $L$ that
have complex multiplication. In other words, if we fix any $L$, the problem of checking when
an elliptic curve over $L$ has CM becomes trivial from a complexity viewpoint: 
pre-compute this list of $j$-invariants for
the field and check if the curve is one of them. The pre-computation cost though prohibitive
is still a computation that requires only $O(1)$ time. For instance, the list for $L = \mathbb{Q}$
is given in \S7.2 of \cite{coh93}. This is why we insist on the field being part of the input.\\

\begin{Rem} The $j$-invariants of elliptic curves with CM are called {\em singular moduli}, and these
enjoy many nice properties. They turn out to be algebraic integers and generate dihedral
extensions of $\mathbb{Q}$. Furthermore, in an important paper Gross-Zagier (\cite{gz85})
derived a formula for the prime ideal factorization of $j(\tau_1) - j(\tau_2)$ where $\tau_1,\tau_2$
generate maximal quadratic orders with coprime discriminants. Such numbers are divisible by
many primes of small norm. There is even a conjectural extension of this work to the case
where the $\tau_i$ do not generate maximal orders; see \cite{hut98}. We utilize some of these
properties in \S\ref{sec_determ:this}.
\end{Rem}

\section{A Direct Approach}\label{sec_direct:this}
We can turn the results of Theorems \ref{thm_cm_crit:this} and \ref{thm_hilb:this} into an algorithm
for checking if an elliptic curve has CM as follows. First compute the Hilbert class
polynomials $H_D = \prod(x - j(\alpha))$, where $\alpha$ runs over the quadratic numbers associated
to the reduced quadratic forms of (negative) discriminant $D$. Next we check if the $j$-invariant of the elliptic
curve is a root of this polynomial. If so, we know that $E$ has CM by an order of discriminant $D$.
This computation can be done in $|D|^{O(1)}$ time (cf. \cite{sch85} \S4). One does this for each 
$D \equiv 0,1 \mod 4$ until the degree of $H_D$ exceeds the degree of the field of definition of the
elliptic curve. At this point we declare that the curve does not have CM.\\

The problem with the above approach is: When do we stop trying new discriminants? The Brauer-Siegel theorem says
that $h(D)$ grows roughly as $|D|^{\frac{1}{2}}$, but this bound is not effective. 
We need an explicit lower bound for the class number
in terms of the discriminant to be able to decide when to stop. This is a hard problem,
first studied by Gauss. Only recently
the following explicit bound was proved by Gross, Zagier, Goldfeld and Osterl{\'e}
(see \cite{zag84,gz86}):

\begin{Thm} If $D$ is a negative fundamental discriminant, then
\begin{align*}
	h(D) > \begin{cases}
	\frac{1}{7000} \ln(|D|) \prod_{p|D}\biggl( 1 - \frac{\lfloor 2\sqrt{p}\rfloor}{p+1}\biggr),
			&\text{ if } \gcd(D,5077) \neq 1 \\
	\frac{1}{55}\ln(|D|)\prod_{p|D}\biggl( 1 - \frac{\lfloor 2\sqrt{p}\rfloor}{p+1}\biggr)
			&\text{ otherwise.}
	\end{cases}
\end{align*}
\end{Thm}

Using the fact that the class number of an order is a multiple of the class number of the quadratic
field associated to it, and the observation that if $D$ has $t$ prime factors then $2^{t-1}~|~h(D)$
(by Gauss's genus theory), we obtain an effective lower bound on $h(D)$. 
This results in a method whose running time is exponential in the degree of the field.

\section{The Randomized Algorithm}\label{sec_rand:this}
The randomized algorithm is based on the observation that if $E/L$ has CM, then there
is an abundance of supersingular primes. This differs from the case where $E$ does not have CM. 
We describe the algorithm first:\\

{\bf Input: } A number field $L$ and $E: Y^2Z = X^3 + AXZ^2 + BZ^3$, with $A,B \in L$.\\
{\bf Steps: }
\begin{enumerate}
\item If $j_E$ is not an algebraic integer, output ``E does not have CM.''
\item Pick a prime $p$ at random in the interval $\mathcal{I} = [2\cdots (h\exp(n^{2+\epsilon})
\max\{\mathbf{w}(A),\mathbf{w}(B)\})^c]$, where $c,h$ and $\epsilon$ are positive
constants and $n = [L:\mathbb{Q}]$.
\item Find the decomposition of $(p) = \prod_{i} \mathfrak{P}_i^{e_i}$, where $\mathfrak{P}_i$
are prime ideals of $\mathscr{O}_L$ (the ring of integers of $L$). If this step fails go back to step (2).
\item Choose a prime in this factorization uniformly at random (say) $\mathfrak{P}$, treating
the $e_i$ copies of $\mathfrak{P}_i$ as distinct.
\item If $N_{L/\mathbb{Q}}\mathfrak{P}$ lies outside the interval $\mathcal{I}$ then go to step (2).
\item With probability $\frac{1}{\deg \mathfrak{P}}$ proceed with the next step; otherwise,
return to step (2).
\item Compute the reduction $\tilde{E}$ of $E ~\mathrm{mod}~ \mathfrak{P}$. If this step does not
suceed return to step (2).
\item Compute $a_{\mathfrak{P}}$, the trace of the Frobenius endomorphism of $\tilde{E}$.
\item If $a_{\mathfrak{P}} = 0 \mod p$ then output ``E probably has CM''; otherwise, 
output ``E probably does not have CM.''
\end{enumerate}

First we argue that all the steps can be done efficiently, and also bound the probability
of failure in some of the steps. Step (1) can be done by computing the minimal polynomial
of $j_E$ and checking if it is monic with integer coefficients. This can be done
in polynomial time \cite{len91}. Step (2) can be done efficiently using our
source of random bits and randomized primality testing methods. To find the splitting
of the prime $p$ we make use of Theorem 4.8.13 in \cite{coh93}, which leads to
a randomized polynomial time algorithm. This algorithm not only provides us with the
prime factorization $(p) = \prod_i \mathfrak{P}_i^{e_i}$ but also gives us
the isomorphism $\mathscr{O}_L/\mathfrak{P} \cong \mathbb{F}_{p^d}$, where $\mathfrak{P} \supseteq (p)$
is a prime and $d = \deg(\mathfrak{P})$. The isomorphism can be used to compute the
reduction of the curve in step (7). The prime decomposition method we suggest will fail if
the prime $p$ divides the index $[\mathscr{O}_L:\mathbb{Z}[\theta]]$, where $\theta = j_E$
(note that $\theta$ is an algebraic integer as a consequence of the check made at step (1)).
The number of primes for which this failure can occur is bounded by the number of
primes that divide the discriminant of the order $\mathbb{Z}[\theta]$. Since
this order has a basis of the form $1,\theta,\theta^2,\cdots,\theta^{n-1}$, its
discriminant is that of its minimal polynomial $T(x) = x^n + a_1x^{n-1} + \cdots +a_n$.
Using the Hadamard bound, we see that the number of primes dividing the discriminant is
bounded by $\log\bigl( ( \sum_i (na_i)^2)^{2n-1}\bigr)$ which is still polynomial
in the input length. The reduction of the elliptic curve can be done in step (7)
if $p \not| N_{L/\mathbb{Q}}\Delta_E$ and this again excludes only a few primes.
Thus, if $c$ and $h$ are large enough the probability that we pick a prime for which either step (3)
or (7) fails will be negligible. Step (8) can be done in polynomial time using, for instance,
Schoof's algorithm \cite{sch85}.\\

We now explain the reason for sampling the primes as we do in steps (2) - (5). We wish
to pick primes $\mathfrak{P}$ uniformly at random from the primes of $\mathscr{O}_L$
whose norm lies in the interval $\mathcal{I}$. The sampling method we use is acceptance-rejection
sampling and this ensures that we pick primes according to our requirement.\\

Firstly, if $E$ has CM then
its $j$-invariant is an algebraic integer (Theorems \ref{thm_cm_crit:this} and \ref{thm_hilb:this}), and step (1) checks that this
holds.
Next, we argue that if $E$ has complex multiplication then with non-negligible probability
the algorithm will output that $E$ probably has CM. 
For this we need a theorem of Deuring
(\cite{lan87} Chapter 13 \S4):

\begin{Thm}[Deuring]\label{thm_deuring_ss:this} 
Let $E/L$ be an elliptic curve with complex multiplication by an
order $\mathscr{O}_E$ of an imaginary quadratic field $K$. Let $\mathfrak{P}$ be
a prime ideal over the rational prime $p$. Assume that $E$ has good reduction at $\mathfrak{P}$.
Then $E \mod \mathfrak{P}$ is supersingular iff $p$ either ramifies or remains inert
in $K$.
\end{Thm}

Let $E$ be an elliptic curve over a finite field $\mathbb{F}_{p^d}$. Then $E$ is
supersingular iff it has no $p$-torsion points. This is equivalent to the trace
of the $p^d$-power Frobenius endomorphism being a multiple of $p$ (\cite{sil86} V. Ex. 5.10).
Thus step (9) checks if $E$ has supersingular reduction at the prime $\mathfrak{P}$.\\

Suppose $E/L$ is a curve with complex multiplication by an order in the imaginary quadratic field
$K = \mathbb{Q}(\sqrt{D})$ (where $D$ is the discriminant of $K$.) Then by Theorem \ref{thm_deuring_ss:this}
the primes where $E$ has supersingular reduction are precisely those primes that are either
ramified or inert in $K$. The primes that ramify are those that divide the discriminant $D$,
and the primes $p$ that remain inert are those for which $\bigl( \frac{D}{p}\bigr)=-1$.
This immediately suggests that the proportion of such primes can be worked out by
choosing primes in certain arithmetic progressions mod $D$. However, since the discriminant
of the field $K$ depends on the input, we need a result that is {\em uniform} in the modulus $D$. Indeed,
using quadratic reciprocity and the uniform prime number theorem for arithmetic progressions (\cite{dav00}
Chapter 20) one can show the following theorem:
\begin{Thm} \label{thm_unif:this}
Define
\begin{align*}
	\pi_0(x) = \sharp\biggl\{ p \leq x ~:~ \biggl(\frac{D}{p}\biggr) = -1\biggr\}
\end{align*}
and let $\delta > 0$ be fixed. Then there is a positive effective constant $c > 0$ depending
on $\delta$ such that if $|D| \leq (\log x)^{1-\delta}$ then
\begin{align*}
	\pi_0(x) = \frac{1}{2}\Li(x) + O\bigl( x e^{-c\sqrt{\log x}}\bigr)
\end{align*}
uniformly in $D$.
\end{Thm}

To apply Theorem \ref{thm_unif:this} we need to ensure that $|D| \leq \log^{1-\delta} x$.
In other words, we need to pick primes in an interval which is longer than $\exp(|D|^{\frac{1}{1-\delta}})$
for some $\delta > 0$. At this point we apply Siegel's theorem to get a bound on $|D|$ in
terms of the degree of the field over which $E$ is defined. We use Siegel's theorem, even
though it is ineffective, because the ineffectiveness affects only the error term in
the success probability of the algorithm. This does not affect the implementation
of the algorithm.

\begin{Thm}[Siegel] For each $\epsilon > 0$ there is a constant (ineffective) $c > 0$ such
that the class number $h(-D)$ satisfies
\begin{align*}
	h(-D) \geq c D^{\frac{1}{2} - \epsilon}.
\end{align*}
\end{Thm}

By Theorem \ref{thm_hilb:this} we have that $[L:\mathbb{Q}] = h(-D)$,
where $-D$ is the discriminant of the order by which $E$ has CM. By Siegel's theorem
we get that $D \leq c'[L:\mathbb{Q}]^{2+\epsilon}$, where $c'$ is a positive constant 
depending on $\epsilon$. Thus picking primes that are at least $\exp(c'[L:\mathbb{Q}]^{2+\epsilon})$
will ensure (Theorem \ref{thm_unif:this}) that we have a positive density of supersingular primes.
In summary, we have proved the following theorem:\\
\begin{Thm}\label{thm_prob_cm:this}
Fix any $\epsilon > 0$ and let $E/L$ be an elliptic curve with CM.
If $p$ is a prime picked uniformly at random in an interval containing 
$[2\cdots \exp([L:\mathbb{Q}]^{2+\epsilon})]$
and $E$ has good reduction at $\mathfrak{P} \supseteq (p)$, then the probability
that $E$ has supersingular reduction at $\mathfrak{P}$ is at least $\frac{1}{2}+o(1)$, the
error term being ineffective.
\end{Thm}

We have shown that about $\frac{1}{2}$ of the {\em rational primes} give us primes of supersingular
reduction for $E$. But our algorithm selects primes $\mathfrak{P}$ of $\mathscr{O}_L$ that are most
likely degree $1$ primes. We need to ensure that this somehow does not bias against the primes of
supersingular reduction for $E$. To argue this we consider the following diagram of fields:
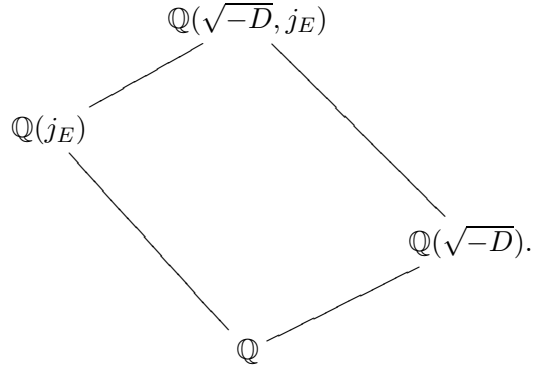
\begin{align*}
\xymatrix{
	& \mathbb{Q}(\sqrt{-D},j_E) & \\
\mathbb{Q}(j_E) \ar@{-}[ur]& & \\
	& & \ar@{-}[uul]\mathbb{Q}(\sqrt{-D}). \\
	& \ar@{-}[uul]\ar@{-}[ur]\mathbb{Q} & 
}
\end{align*}

All extensions in the diagram are galois, except possibly the extension $\mathbb{Q}(j_E)/\mathbb{Q}$
(\cite{shi71} Theorem 5.7). Now since $\mathbb{Q}(\sqrt{-D},j_E)/\mathbb{Q}(j_E)$
is a degree $2$ extension, the Chebotarev density theorem tells us that 
\begin{align}\label{eqn_dens_in:this}
\sharp\{\mathfrak{P}~:~
N_{\mathbb{Q}(j_E)/\mathbb{Q}}\mathfrak{P} \leq x, \deg \mathfrak{P} = 1 
\text{ and }\mathfrak{P}\text{ remains inert in } \mathbb{Q}(\sqrt{-D},j_E)\} \sim \frac{1}{2}\Li(x).
\end{align}
If $\mathfrak{P}$ is a (degree $1$) prime of $\mathbb{Q}(j_E)$ that remains inert in $\mathbb{Q}(\sqrt{-D},j_E)$
then its norm (a rational prime) remains inert in $\mathbb{Q}(\sqrt{-D})$. Such a prime $\mathfrak{P}$ 
is a supersingular prime if $E$ has good reduction at $\mathfrak{P}$. Thus
we have shown that half of the degree $1$ primes of $L$ are indeed primes of supersingular
reduction for $E$. In particular, if our algorithm is given an elliptic curve with CM, then it
outputs ``E probably has CM'' with probability $\geq \frac{1}{2} + o(1)$.\\

Now suppose $E/L$ does {\em not} have CM. Then we show that the probability that we pick
a prime $p$, where $E$ has supersingular good reduction at a prime above $p$ goes to $0$.
For this we use a result of Serre (\cite{ser81} \S8) that says:
\begin{Thm}\label{thm_ser_noncm:this} 
Let $E/L$ be an elliptic curve that does not have CM and let
\begin{align*}
\pi_{E,0} = \sharp\{\mathfrak{P}~:~\mathfrak{P}\text{ a prime of }\mathscr{O}_L,
 N_{L/\mathbb{Q}}\mathfrak{P} \leq x, E 
\text{ has supersingular reduction at } \mathfrak{P}\}.
\end{align*}
Then for $\delta > 0$
\begin{align*}
	\pi_{E,0} = O\biggl( \frac{x}{(\log x)^{\frac{3}{2}-\delta}}\biggr).
\end{align*}
The implicit constant depends only on $\delta$.
\end{Thm}

\begin{Rem} Serre states his theorem only for elliptic curves over $\mathbb{Q}$
but the proof works for elliptic curves over number fields too. We sketch a proof of
a weaker form of Theorem \ref{thm_ser_noncm:this} in \S\ref{sec_determ:this}.
There are stronger versions of this result, most notably due to Noam Elkies with some restrictions
on the number field \cite{elk91}, but the weaker version is sufficient for our purpose.
For curves defined over $\mathbb{Q}$, 
a famous conjecture of Lang and Trotter predicts that $\pi_{E,0} \sim C_E \frac{\sqrt{x}}{\log x}$
where $C_E$ is a constant depending on $E$ (\cite{ltr76}).
\end{Rem}

Theorem \ref{thm_ser_noncm:this} immediately gives us the following result:
\begin{Thm} Suppose $E/L$ is an elliptic curve that does not have CM. If $\mathfrak{P}$
is a prime picked uniformly at random among those whose norm lies in the interval $[2\cdots x]$,
then the probability that $E$ has supersingular reduction at $\mathfrak{P}$
tends to $0$ with $x$.
\end{Thm}

Putting Theorems \ref{thm_prob_cm:this}, \ref{thm_ser_noncm:this} and
the remarks following Theorem \ref{thm_prob_cm:this} together, we
see that if $E/L$ has CM then the output of the algorithm is correct with probability $\frac{1}{2}+o(1)$,
and if $E/L$ does not have CM then the output is correct with probability $1 - o(1)$.
This shows that we have a two-sided error randomized polynomial time algorithm
for checking when an elliptic curve over a number field has CM.
If one needs to improve the confidence of the algorithm, then one can use the
standard boosting idea of repeating the algorithm independently many times and taking the
majority vote (cf. \cite{pap95} Corollary to Lemma 11.9).\\

In the Appendix we tabulate the ratio of supersingular primes to all the primes, considering
only primes of norm $\leq 10^5$, for certain curves. 
One sees that for curves with CM, this
ratio is already close to $\frac{1}{2}$, and for curves without CM it is very small.

\subsection{Finding the discriminant of $\mathbf{End(E)}$} Suppose $E/L$ is an
elliptic curve with CM. Then even at the primes where $E$ has non-supersingular
good reduction, the trace of Frobenius gives important information. The following
theorem of Deuring is the main tool we use (\cite{lan87} Chapter 13 \S4, Theorem 12):

\begin{Thm}[Deuring] \label{thm_cm_ord:this}
Let $E/L$ be an elliptic curve with CM by $\mathscr{O}_E$, an
order in an imaginary quadratic field $K$. Assume that $p$ is a rational prime
that splits completely in $K$ and that $\mathfrak{P} \supseteq (p)$ is a prime of $L$
above $p$. Suppose that $E$ has good non-supersingular reduction $\tilde{E}$
at $\mathfrak{P}$ and that $p$ does not divide the index $[\mathscr{O}_K:\mathscr{O}_E]$
($\mathscr{O}_K$ is the ring of integers of $K$). Then $\mathrm{End}(E) \cong \mathrm{End}(\tilde{E})$.
\end{Thm}

Let $E/L$ be a curve with CM by $\mathscr{O}_E$. Suppose we pick
a prime of good reduction $\mathfrak{P}$ of $L$ and find that $a_{\mathfrak{P}} \neq 0 \mod p$
for the reduction $\tilde{E}$ (where $a_{\mathfrak{P}}$ is the trace of Frobenius
on $\tilde{E}$). Then assuming $p$ does not divide the index of $\mathscr{O}_E$
(which happens with high probability), we get from Theorem \ref{thm_cm_ord:this}
that $\mathscr{O}_E = \End(E) \cong \End(\tilde{E})$. 
Since $\tilde{E}$ is an elliptic curve over a finite field $\mathbb{F}_{p^d}$,
($d = $ degree of $\mathfrak{P}$) the $p^d$-power Frobenius endomorphism $\phi$
satisfies
\begin{align}\label{eqn_char_pol:this}
\phi^2 - a_{\mathfrak{P}} \phi + p^d = 0
\end{align}
as an element of $\End(\tilde{E})$. Since the latter is an order with discriminant $D_{\mathscr{O}_E}$
(say) equation (\ref{eqn_char_pol:this}) implies that
\begin{align}\label{eqn_sqrmul:this}
a_{\mathfrak{P}}^2 - 4p^d = m_{\mathfrak{P}}^2 D_{\mathscr{O}_E}
\end{align}
for some $m_{\mathfrak{P}} \in \mathbb{Z}.$
Since $\tilde{E}$ is not supersingular this quantity is never $0$. The idea is to pick
different primes $\mathfrak{P}_i$ (assume that the reduction of the curve is non-supersingular),
and compute the quantities $w_i = a_{\mathfrak{P}_i}^2 - 4p^d$ and $\gcd(w_i)$. We hope this 
gives us $D_{\mathscr{O}_E}$. However, we do not know how to argue that
the $\gcd(w_i)$ quickly converge to the discriminant. In experiments, two trials
were sufficient in every case we tested. Another piece of information that equation (\ref{eqn_sqrmul:this})
and Hasse's bound yield is this. If $4p^d < |D_{\mathscr{O}_E}|$, then the hypotheses
of Theorem \ref{thm_cm_ord:this} must fail. Thus the curve either
has bad reduction, or supersingular reduction, or $p$ must divide the index
of the order $\mathscr{O}_E$. In the last case it turns out
that the endomorphism ring of $\tilde{E}$ is an order of index $[\mathscr{O}_K:\mathscr{O}_E]/p^r$, 
where $p^r$ is the largest power of $p$ dividing the index of $\mathscr{O}_E$. Thus we
get some information about the index of $\mathscr{O}_E$.
If, on the other hand, $E$ does not have CM, then the $w_i$ should behave randomly
and we should get $\gcd(w_i) = 1$ very quickly. Again, we are unable to prove this.

\begin{Rem}\label{rem_oneside:this} We can use the ideas here to  make the error in the 
randomized algorithm one-sided. Taking a bunch of primes $\mathfrak{P}_i$ and reducing the curve
we can find the quantity $w_i$ (for those primes of ordinary reduction). 
If $\gcd(w_i) = 1$, then we know for certain
that the curve does not have CM. However, we cannot {\em prove} that if $E$ does
not have CM, then this will happen for a reasonable number of primes $\mathfrak{P}_i$.
The method in \cite{cnst98} also incorporates a similar idea, but in their
proof (of Theorem 3) they claim, in essence, that the $w_i$ behave like random
numbers without proof. Our algorithm in \S\ref{sec_rand:this} has two-sided error,
but its behavior is rigorously proved. If one uses the one-sided error version, then
its running time analysis needs the heuristic assumption that the $w_i$ behave like
random numbers if $E$ does not have CM.
\end{Rem}

\section{The Deterministic Algorithm}\label{sec_determ:this}

This method uses the galois representations that
are afforded by the elliptic curve. We briefly describe such galois representations in the next subsection.

\subsection{Galois Representations from Elliptic curves} 
For more on this subject the reader should consult Serre (\cite{ser89}) and also
Silverman (\cite{sil86} III \S7). Let $E/L$ be an elliptic curve and let $\ell$ be a prime.
The set of $\ell$-torsion points on $E$ is
\begin{align*}
	E[\ell] = \{P \in E(\mathbb{C}) ~:~ \ell P = \infty\},
\end{align*}
where $\infty$ is the identity on $E$. It is known that 
$E[\ell] \cong (\mathbb{Z}/\ell\mathbb{Z})\times(\mathbb{Z}/\ell\mathbb{Z})$ (\cite{sil86} III \S6.4).
Let $G_{L} = \Gal(\overline{L}/L)$ be the absolute galois group of $L$. If $K \supseteq L$ is
a galois extension, then $G_L$ acts on $E(K)$ (the points on $E(\mathbb{C})$ with coordinates in $K$) by 
sending the point $(x:y:z)$ to $(x^{\sigma}:y^{\sigma}:z^{\sigma})$ for $\sigma \in G_L$.\\

$G_L$ also acts on $E[\ell]$ since the multiplication by $\ell$ maps are defined over $L$. Thus we
get a map
\begin{align*}
\rho_{\ell} : G_L \rightarrow \Aut(E[\ell]) \cong \GL_2(\mathbb{F}_{\ell}).
\end{align*}
This is a continuous group homomorphism (with profinite topology on $G_L$ and discrete
topology on $\GL_2(\mathbb{F}_{\ell})$) and gives us a representation of $G_L$.
Now if $\sigma \in \Gal(\overline{L}/L(E[\ell]))$ then it acts trivially on $E[\ell]$.
Thus the representation factors through the extension $L(E[\ell])$ and we get a representation
of $\Gal(L(E[\ell])/L)$:
\begin{align*}
\rho_{\ell}: \Gal(L(E[\ell])/L) \rightarrow \GL_2(\mathbb{F}_{\ell}).
\end{align*}
The representation is clearly injective. It turns out that $\mathrm{Im} ~\rho_{\ell}$
depends critically on whether $E$ has CM or not. We discuss this next.

\subsection{Image of $\rho_{\ell}$ if $E$ does not have CM} Suppose $E/L$ does not have
CM. Then a famous theorem of Serre (\cite{ser72}) says the following:

\begin{Thm}\label{thm_im_noncm:this}
Let $E/L$ be an elliptic curve that does not have CM. Then for all large enough
primes $\ell$, the representation $\rho_{\ell}$ is surjective, i.e., $\rho_{\ell}(G_L) = \GL_2(\mathbb{F}_{\ell})$.
This means that
\begin{align*}
	\Gal(L(E[\ell])/L) \cong \GL_2(\mathbb{F}_{\ell})
\end{align*}
for all but finitely many primes $\ell$.
\end{Thm}

We illustrate the power of this theorem by sketching a proof of the following result.

\begin{Cor} Let $E/L$ be an elliptic curve without complex multiplication. Then
\begin{align*}
\sharp\{\mathfrak{P}~:~\mathfrak{P}\text{ a prime of } \mathscr{O}_L, 
N_{L/\mathbb{Q}}\mathfrak{P} \leq x, E ~\mathrm{mod}~ \mathfrak{P}
\text{ is supersingular}\} = o(\Li(x)).
\end{align*}
\end{Cor}
\begin{Proof} Fix a prime $\ell$. 
We need the following fundamental compatibility between the
Frobenius at a prime $\mathfrak{P}$ of $L$ and the Frobenius on $E ~\mathrm{mod}~ \mathfrak{P}$ via
the representation $\rho_{\ell}$. Suppose $\mathfrak{P}$ is a prime where $E$ has good reduction,
and assume that $\mathfrak{P}$ does not divide the discriminant of $L$. Then
\begin{align*}
\Tr(\rho_{\ell}(\mathrm{Frob}_{\mathfrak{P}})) \equiv a_{\mathfrak{P}} \mod \ell,
\end{align*}
where $a_{\mathfrak{P}}$ is the trace of Frobenius on the curve.\\

Let $\ell_0$ be such that for all primes $\ell \geq \ell_0$ the representation $\rho_{\ell}$
coming from $E$ is surjective. Now for any prime $\ell \geq \ell_0$  we have that
$\Gal(L(E[\ell])/L) \cong \GL_2(\mathbb{F}_{\ell})$. Let $S_0$ be the set of primes 
\begin{align*}
\{\mathfrak{P} ~:~ E \text{ has good reduction at }\mathfrak{P}\text{ and } a_{\mathfrak{P}} = 0\}.
\end{align*}
Note that the set $S_0$ contains all the degree $1$ primes where $E$ has supersingular reduction.
The Chebotarev density theorem says that the density of primes $\mathfrak{P}$ such
that $\Tr(\rho_{\ell}(\mathrm{Frob}_{\mathfrak{P}})) \equiv 0 \mod \ell$ is exactly the ratio
\begin{align*}
r_{\ell} = \frac{\sharp\{\text{Trace }0\text{ conjugacy class of }\GL_2(\mathbb{F}_{\ell})\}}
	{\sharp\GL_2(\mathbb{F}_{\ell})}.
\end{align*}
A quick calculation shows that $r_{\ell} \ll \frac{1}{\ell}$. Now
\begin{align*}
	\lim_{\ell \rightarrow \infty} r_{\ell} = 0,
\end{align*}
proving that the density of the set $S_0$ is $0$ (counted by norm). The set
of primes of $L$ which are of degree $> 1$ are already density $0$, when we are counting by norm.
So that even among the degree $1$ primes there is only a density $0$ subset
where $E$ has supersingular reduction.
\end{Proof}
\subsection{Image of $\rho_{\ell}$ if $E$ has CM} If $E/L$ has CM we have, from
the theory of complex multiplication (\cite{sil94} Chapter II Theorem 2.3), the
following result.

\begin{Thm}\label{thm_cm:this} Let $E/L$ be an elliptic curve that
has complex multiplication by an order $\mathscr{O}_E$ in $\mathbb{Q}(\sqrt{D})$ ($D < 0$)
and let $\ell$ be a prime. Then $L(\sqrt{D},E[\ell])/L(\sqrt{D})$ is an abelian extension.
\end{Thm}

Now consider the following diagram of fields:
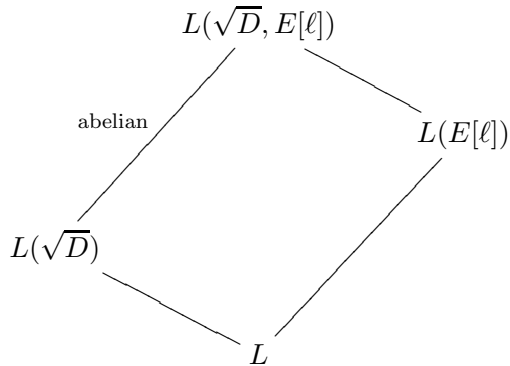
\begin{align*}
\xymatrix{
 & L(\sqrt{D}, E[\ell])& \\
 & & L(E[\ell]) \ar@{-}[ul]\\
L(\sqrt{D}) \ar@{-}[uur]^{\text{abelian}} & & \\
 & L \ar@{-}[ul] \ar@{-}[uur] &
}
\end{align*}
The group $\Gal(L(\sqrt{D},E[\ell])/L(\sqrt{D}))$ is an abelian subgroup
of $\Gal(L(\sqrt{D},E[\ell])/L)$, furthermore, it has index $2$. This implies
that $\Gal(L(\sqrt{D},E[\ell])/L)$ is solvable. Therefore $\Gal(L(E[\ell])/L)$,
being a quotient of a solvable group, is also solvable. We have thus proved:
\begin{Thm}\label{thm_im_cm:this}
Suppose $E/L$ is an elliptic curve with complex multiplication, and $\ell$ a prime.
Then $\mathrm{Im}~\rho_{\ell}$ is solvable.
\end{Thm}

\subsection{The algorithm} The idea is to use Theorems \ref{thm_im_noncm:this} and
\ref{thm_im_cm:this} to check if $E$ has CM. We pick $\ell \geq 5$
and large enough so that if $E$ did not have CM then $\rho_{\ell}$ would have
to be surjective. Since $\SL_2(\mathbb{F}_{\ell})$, a subgroup of $\GL_2(\mathbb{F}_{\ell})$,
is not solvable for $\ell \geq 5$, $\GL_2(\mathbb{F}_{\ell})$ is not solvable for $\ell \geq 5$.
In summary, if $\ell$ is large enough, then $\Gal(L(E[\ell])/L)$ is solvable
iff $E/L$ has complex multiplication. The extension $L(E[\ell])/L$ 
is of degree $\leq \sharp\GL_2(\mathbb{F}_{\ell}) = (\ell^2-1)(\ell^2-\ell)$.
Solvability of this extension can be checked in polynomial time, provided,
$\ell$ is bounded polynomially in the input length. This can be done by
computing the $\ell$ division polynomial of $E$ and using the
algorithm of Landau and Miller \cite{len91}. To complete the description
of the algorithm we need to decide how large an $\ell$ to take. The following theorem
of Masser and W{\"u}stholz \cite{mwu93} allows us to do that.

\begin{Thm} There are absolute constants $c,\gamma$ ($\gamma$ is effectively computable)
with the following properties. Suppose $E$ is an elliptic curve of Weil height $h$ defined
over a number field $L$ of degree $d$, and assume that $E$ does not have complex multiplication.
\begin{enumerate}
\item If $\ell > c(\max\{d,h\})^{\gamma}$, then $\rho_{\ell}(G_L)$ contains the special linear
group $\SL_2(\mathbb{F}_{\ell})$.
\item If, further, $\ell$ does not divide the discriminant of $L$, then $\rho_{\ell}(G_L) 
= \GL_2(\mathbb{F}_{\ell})$.
\end{enumerate}
\end{Thm}

If $\rho_{\ell}$ contains $\SL_2(\mathbb{F}_{\ell})$ for $\ell \geq 5$ then it is already
non-solvable, thus we get the following result:

\begin{Thm}\label{thm_determ:this}
There are absolute constants $c,\gamma$ ($\gamma$ effective) with the
following property. 
Suppose $E/L$ is
an elliptic curve of Weil height $h$, $d = [L:\mathbb{Q}]$, and
$\ell > \max\{c(\max\{d,h\})^{\gamma},5\}$ is a prime. Then $E$ has complex multiplication
iff $\Gal(L(E[\ell])/L)$ is solvable.
\end{Thm}

Since the Weil-height of the elliptic curve is bounded polynomially by the input length,
we get a deterministic polynomial time algorithm to test if $E/L$ has complex multiplication.
Unfortunately, the constant in the running time has not yet been made effective. Serre has
conjectured that the lower bound on the primes for which $\rho_{\ell}$ is
surjective for curves without CM over $L$ should only depend on $L$ and not
on the curve \cite{ser72} \S4.3. For all the curves (without CM) we tested $\ell = 5$
or $7$ already gave non-solvable extensions. It must be noted however, that there are
curves over $\mathbb{Q}$ for which $\rho_{\ell}$ is not surjective if $\ell < 47$.\\

{\bf Acknowledgements: } I would like to thank Eric Bach, Nigel Boston, Rohit Chatterjee, Ken Ono and
Gisbert W{\"ustholz} for extremely useful discussions and suggestions. I am especially grateful to 
Nigel for suggesting to look at the image of Galois and to Eric for help with the acceptance-rejection
sampling method.

\pagebreak

\appendix
\section*{Appendix}

In this appendix we tabulate the ratios of supersingular
to ordinary primes for some elliptic curves. In each case
if $E/L$ is an elliptic curve, we computed the ratio
\begin{align*}
\frac{\pi_{E,0}(10^5)}{\pi(10^5)} &=
\frac{\sharp\{\mathfrak{P} \text{ prime of } \mathscr{O}_L~:~ 
	\mathfrak{P}\text{ relatively prime to }\Delta_E\text{ and } N_{L/\mathbb{Q}}(\mathfrak{P}) \leq 10^5\}}
{\sharp\{\mathfrak{P} \text{ prime of }\mathscr{O}_L~:~N_{L/\mathbb{Q}}\mathfrak{P} \leq 10^5\}}.
\end{align*}
All our computation was done using MAGMA version 2.10 \cite{bc03}.\\

In Table \ref{tbl_cm:this} we give
the results for elliptic curves with complex multiplication. To prepare this
table we picked elliptic curves with CM by the maximal orders of $\mathbb{Q}(\sqrt{-p})$
with $p$ a prime in the range $50 \leq p \leq 100$. We ignored those $p$ for which the
class number of $\mathbb{Q}(\sqrt{-p})$ is $1$, since these curves are then defined over $\mathbb{Q}$.
The entries in the table are listed in increasing order of the prime $p$.\\

\begin{table}
\caption{Proportion of Supersingular primes for CM curves}
\label{tbl_cm:this}
\begin{tabular}{|c|c|c|} \hline
Discriminant $D$ of $\End(E)$ & Degree of Number field $L$ & $\frac{\pi_{E,0}(10^5)}{\pi(10^5)}$ \\\hline
$-4\times 53$ & 6  & 0.5043 \\\hline
$-59$  & 3  & 0.5073 \\\hline
$-4\times 61$ & 6  & 0.5079 \\\hline
$-71$  & 7  & 0.5113 \\\hline
$-4\times 73$ & 4  & 0.5110 \\\hline
$-79$  & 5  & 0.5107 \\\hline
$-83$  & 3  & 0.5088 \\\hline
$-4 \times 89$ & 12 & 0.5234 \\\hline
$-4 \times 97$ & 4  & 0.5040 \\\hline
\end{tabular}
\end{table}

In Table \ref{tbl_noncm:this} we give the results for elliptic curves without complex multiplication over
a degree $5$ number field.
The table was prepared by picking random monic polynomials of degree $5$ and using
a root of the polynomial as the $j$-invariant of the elliptic curve. We verified
that these curves do not have CM by using the criterion described in Remark \ref{rem_oneside:this}.
We see that the results of these experiments are consistent with Theorems 
\ref{thm_prob_cm:this} and \ref{thm_ser_noncm:this}.

\begin{table}
\caption{Proportion of Supersingular primes for Non-CM curves}
\label{tbl_noncm:this}
\begin{tabular}{|c|c|} \hline
Minimal polynomial of $j$-invariant & $\frac{\pi_{E,0}(10^5)}{\pi(10^5)}$ \\\hline
$x^5 - 12x^4 - 65x^3 - 33x^2 -22x -51$   & 0.0032 \\\hline
$x^5 - 78x^4 + 28x^3 + 14x^2 - 92x + 19$ & 0.0036 \\\hline
$x^5 + 25x^4 + 7x^3 + 25x^2 + 96x + 92$  & 0.0035 \\\hline
$x^5 + 71x^4 - 71x^3 + 41x^2 + 61x + 93$ & 0.0034 \\\hline
$x^5 + 23x^4 + 84x^3 - 17x^2 - 36x + 62$ & 0.0031 \\\hline
$x^5 - 94x^4 - 74x^3 + 78x^2 + 51x - 10$ & 0.0033 \\\hline
$x^5 + 79x^4 + 97x^3 + 5x^2 - 78x - 39$  & 0.0033 \\\hline
$x^5 + 68x^4 - 17x^3 + 99x^2 - 34x - 93$ & 0.0025 \\\hline
\end{tabular}
\end{table}

\end{document}

%% file: MTemplate.tex
\usepackage{fullpage}
\usepackage{amssymb}
\usepackage{amsmath}
\usepackage{amsfonts}
\usepackage{amsxtra}
\usepackage{amscd}
\usepackage{amsthm}
\usepackage{epsfig}
\usepackage{enumerate}
\usepackage{float}
\usepackage{xy}

\xyoption{all}

\theoremstyle{plain} \newtheorem{Thm}{Theorem}[section]
\theoremstyle{plain} \newtheorem{Cor}[Thm]{Corollary}
\theoremstyle{plain} 
\theoremstyle{plain} 
\theoremstyle{definition} \newtheorem{Def}[Thm]{Definition}
\theoremstyle{definition} 
\theoremstyle{remark} \newtheorem{Rem}[Thm]{Remark}
\newenvironment{Proof}{{\bf Proof :}}{$\Box$\\}

\newcommand{\Gal}{\mbox{$\mathrm{Gal}$}}
\newcommand{\Tr}{\mbox{$\mathrm{Tr}~$}}

\newcommand{\SL}{\mbox{$\mathrm{SL}$}}
\newcommand{\GL}{\mbox{$\mathrm{GL}$}}

\newcommand{\Li}{\mbox{$\mathrm{Li}$}}

\newcommand{\Aut}{\mbox{$\mathrm{Aut}$}}
\newcommand{\End}{\mbox{$\mathrm{End}$}}

\newcommand{\newfontobj}[2]{
  \newcommand{#1}[1]{
    \expandafter\def\csname##1\endcsname{{#2 ##1}}}}

\newfontobj{\class}{\rm}
\newfontobj{\lang}{\bf}
\newfontobj{\oper}{\rm}

\class{PSPACE}		
\class{AM}		
\class{MA}
\class{NP}
\class{UP}
\class{P}
\class{RP}
\class{BPP}
\class{DTIME}
\class{ZPP}
\class{EXPSPACE}
\class{coNP}
\class{coRP}
\class{coAM}
\class{PH}
\class{co}
\class{BP}	

\lang{HN} 		
\lang{SAT}
\lang{PIM}
\lang{CSG}
\lang{ESC}
\lang{IP}
\lang{PRIMES}

\oper{Pr}    
\oper{E}     
\oper{Ord}   
\oper{li}

\floatstyle{ruled}
\newfloat{Algorithm}{thp}{loa}[section]


%% file: ComplexMult.bbl
\begin{thebibliography}{cdx04}
\bibitem[AtMor93]{atmor93} Atkin, A., O., L.; Morain, F.; {\em Elliptic curves and primality proving},
Math. Comp., {\bf 61}, no. 203, 29-68, 1993.

\bibitem[BC03]{bc03}
Bosma, W.; Cannon, J.; 
{\em Handbook of MAGMA functions},
Sydney, 2003.

\bibitem[CNST98]{cnst98}
Chao, J.; Nakamura, O.; Sobataka, K.; Tsujii, S.;
{\em Construction of secure elliptic cryptosystems using CM tests and liftings},
Advances in Cryptology, ASIACRYPT'98 (Beijing), Lecture Notes in Computer Science, 1514,
Springer-Verlag, Berlin, 1998.

\bibitem[Coh93]{coh93}
Cohen, Henri; {\em A course in Computational Algebraic Number Theory}, Graduate Texts
in Math., Vol. 138, Springer-Verlag, 1993.

\bibitem[Dav00]{dav00}
Davenport, Harold; {\em Multiplicative Number Theory}, 3rd ed., revised by Hugh L. Montgomery,
Graduate Texts in Math., vol. 74, Springer-Verlag, 2000.

\bibitem[Elk91]{elk91}
Elkies, Noam, D.; {\em Distribution of Supersingular primes},
Ast{\'e}risque, {\bf 198-200}, 127-132, 1991.

\bibitem[Fel82]{fel82}
Feldman, N., I.; {\em The seventh Hilbert's problem}, Moscow, Moscow State University, 1982.

\bibitem[GZ85]{gz85}
Gross, B.; Zagier, D.; {\em On singular moduli},
J. Reine Angew. Math., {\bf 355}, 191-220, 1985.

\bibitem[GZ86]{gz86}
Gross, B.; Zagier, D.;
{\em Heegner points and derivatives of L-series}, Invent. Math.,
{\bf 84}, no. 2, 225-320, 1986.

\bibitem[Hut98]{hut98}
Hutchinson, Tim; {\em A conjectural extension of the Gross-Zagier formula on singular moduli},
Tokyo J. Math., {\bf 21}, no. 1, 255-265, 1998.

\bibitem[Lan87]{lan87}
Lang, Serge; {\em Elliptic Functions}, 2nd ed., Graduate Texts in Math., vol. 112, Springer-Verlag, 1987.

\bibitem[LTr76]{ltr76}
Lang, Serge; Trotter, Hale, F.;
{\em Frobenius distributions in $\GL_2$-extensions},
Lecture Notes in Math., {\bf 504}, Springer-Verlang, 1976.

\bibitem[Len91]{len91}
Lenstra, Hendrik, W., Jr.;
{\em Algorithms in Algebraic Number Theory},
Bull. Amer. Math. Soc., vol. {\bf 26}, no. 2, 211-244, 1991.

\bibitem[MW{\"u}93]{mwu93}
Masser, D., W.; W{\"u}stholz, G.;
{\em Galois properties of division fields of elliptic curves},
Bull. Lond. Math. Soc., {\bf 25}, 247-254, 1993.

\bibitem[Pap95]{pap95}
Papadimitriou, Christos;
{\em Computational Complexity}, Addison-Wesley, 1995.

\bibitem[Sch85]{sch85}
Schoof, Ren{\'e}; {\em Elliptic curves over finite fields and Computation of square
roots $\mathrm{mod}~ p$}, Math. Comp., vol {\bf 44}, no. 170, 483-494, 1985.

\bibitem[Ser72]{ser72}
Serre, Jean-Pierre;
{\em Propri{\'e}t{\'e}s galoisiennes des points d'ordre fini des courbes elliptiques},
Invent. Math., {\bf 16}, 259-331, 1972.

\bibitem[Ser81]{ser81}
Serre, Jean-Pierre;
{\em Quelques applications du th{\'e}or{\`e}me de densit{\'e} de Chebotarev},
Publ. Math. I.H.E.S, {\bf 54}, 123-201, 1981.

\bibitem[Ser89]{ser89}
Serre, Jean-Piere;
{\em Abelian $\ell$-adic representations and elliptic curves},
with the collaboration of Willem Kuyk and John Labute, 2nd ed., Advanced Book
Classics, Addison-Wesley, 1989.


\bibitem[Shi71]{shi71}
Shimura, Goro; {\em Introduction to the Arithmetic Theory of Automorphic functions},
Iwanami Shoten and Princeton University Press, 1971.

\bibitem[Sil86]{sil86}
Silverman, Joseph; {\em The Arithmetic of Elliptic Curves}, Graduate Texts in Math. Vol. 106,
Springer-Verlag, 1986.

\bibitem[Sil94]{sil94}
Silverman, Joseph; {\em Advanced Topics in the Arithmetic of Elliptic Curves}, Graduate Texts in Math. Vol. 151,
Springe-Verlag, 1994.

\bibitem[Zag84]{zag84}
Zagier, Don, B.;
{\em L-series of elliptic curves, the Birch-Swinnerton-Dyer conjecture, and the class number
problem of Gauss}, Notices Amer. Math. Soc., {\bf 31}, no. 7, 739-743, 1984.
\end{thebibliography}
